\documentstyle[12pt]{article}
\normalbaselines
\begin{document}

\title{About rational-trigonometric deformation}

\vskip 15pt
\author{V.N. Tolstoy\thanks{e-mail: tolstoy@nucl-th.sinp.msu.ru}
}
\date{}
\maketitle
\begin{center}
\vskip -30pt
{\large Institute of Nuclear Physics, Moscow State University\\[3pt]
119992 Moscow \& Russia 
}
\end{center}


\begin{abstract}
We consider a rational-trigonometric deformation in context of
rational and trigonometric deformations. The simplest examples of
these deformations are presented in different fields of
mathematics. Rational-trigonometric differential
Knizhnik-Zamolodchikov and dynamical equations are introduced.
\end{abstract}

\section{Introduction}
In classical mathematics there are three classes of meromorphic
functions: rational, trigonometric, elliptic. According to these
classes there are three types of quantum deformations:
\textit{rational, trigonometric, elliptic}. It turns out we can
also introduce a \textit{rational-trigonometric deformation}. All
these deformations will be called standard.

In this short and sketched paper we discussed the simplest
examples of the standard deformations in the arithmetic (number
theory), geometry, differential calculus, theory functions and Lie
algebras. We also introduce rational-trigonometric differential
Knizhnik-Zamolodchikov (KZ) and dynamical (DD) equations. The
rational-trigonometric differential KZ equations are connected
with a rational-trigonometric classical $r$-matrix \cite{tolsT1}
which is a sum the simplest rational and trigonometric
$r$-matrices depending on spectral parameter. It turns out that
the rational-trigonometric differential KZ (DD) equations are sums
of the rational and trigonometric differential KZ (DD) equations.

\section{Trigonometric, rational and rational-trigonometric deformations}
\label{typing section}
\subsection{Standard deformations of numbers}
Let $z$ be a complex number, $z\in \;{\rm I\!\!\!C}$.

a). The trigonometric deformation (or $q$-deformation) of $z$:
\begin{equation}
z\rightarrow(z)_{q}=\frac{1-q^z}{1-q}~, \qquad
(1)_q=1,\qquad(0)_q=0~, \label{tolsDF1}
\end{equation}
where $q$ is a deformation parameter, $(z)_{q=1}=z$. This
deformation is well-known.

b). The rational deformation (or $\eta$-deformation) of $z$
\cite{tolsT2}:
\begin{equation}
z\rightarrow(z)_{\eta}=\frac{z}{1+\eta(z-1)}~, \qquad
(1)_{\eta}=1,\qquad(0)_{\eta}=0~, \label{tolsDF2}
\end{equation}
where $\eta$ is a deformation parameter, $(z)_{\eta=0}=z$.

c). The rational-trigonometric deformation (or
$(q,\eta)$-deformation) of $z$ \cite{tolsT2}:
\begin{equation}
z\rightarrow(z)_{q\eta}=\frac{(z)_q}{1+\eta(z-1)_q}~, \qquad
(1)_{q\eta}=1,\qquad(0)_{q\eta}=0~, \label{tolsDF3}
\end{equation}
where $q$ and $\eta$ are deformation parameters,
$(z)_{q=1,\eta}=(z)_{\eta}$, $(z)_{q,\eta=0}=(z)_q$.

\textit{Remark}. There is an elliptic deformation of $z$ (for
example, see \cite{tolsSp})

\subsection{Standard deformations of two-dimensional plane}
Let $x$ and $y$ be two commuting coordinate variables, i.e.
\begin{equation}
xy-yx=0\qquad\quad{\rm(the}\;\;(x,y)-{\rm plane)}. \label{tolsDF4}
\end{equation}

a). The trigonometric deformation of the $(x,y)$-plane:
\begin{equation}
xy-qyx=0\qquad\quad{\rm (the\;Manin's\;plane)}.\label{tolsDF5}
\end{equation}

b). The rational deformation of the $(x,y)$-plane:
\begin{equation}
xy-yx=\eta y^2~. \label{tolsDF6}
\end{equation}

c). The rational-trigonometric deformation of the $(x,y)$-plane:
\begin{equation}
xy-qyx=\eta y^2~.\label{tolsDF7}
\end{equation}
All these deformation of the $(x,y)$-plane are well-known.

\textit{Remark}. No elliptic deformation of the $(x,y)$-plane is
known. It is an open problem.

\subsection{Standard deformations of the differential calculus}

Let $\partial_x:=\frac{\partial}{\partial x}$ be the usual
derivative.

a). The trigonometric deformation of $\partial_x$:
\begin{equation}
\partial_x^{(q)}f(x)_{q}=\frac{f(qx)-f(x)}{qx-x}=\frac{f(qx)-f(x)}{x(q-1)}~.
\label{tolsDF8}
\end{equation}

b). The rational deformation of $\partial_x$:
\begin{equation}
\partial_x^{(\eta)}f(x)_{q}=\frac{f(x+\eta)-f(x)}{x+\eta-x}=
\frac{f(x+\eta)-f(x)}{\eta}~.\label{tolsDF9}
\end{equation}
The deformations a) and b) for $\partial_x$ are well-known.

c). The rational-trigonometric deformation of $\partial_x$:
\begin{equation}
\partial_x^{(q\eta)}f(x)_{q}=\frac{f(qx+\eta)-f(x)}{qx+\eta-x}=
\frac{f(qx+\eta)-f(x)}{x(q-1)+\eta}~. \label{tolsDF10}
\end{equation}
All these deformations of $\partial_x$ have the following
properties:

(i) they are defined by the general formulas
\begin{equation}
\partial_x^{(df)}f(x)=\frac{\triangle^{(df)}f(x)}{\triangle^{(df)}x}=
\frac{f(x')-f(x)}{x'-x}~, \label{tolsDF11}
\end{equation}

(ii) they satisfy the Leibniz rule
\begin{equation}
\partial_x^{(df)}f(x)\varphi(x)=\left(\partial_x^{(df)}f(x)\right)\varphi(x)+
f(x')\left(\partial_x^{(df)}\varphi(x)\right)~, \label{tolsDF12}
\end{equation}
where $x'=qx$ for the $q$-deformation, $x'=x+\eta$ for the
$\eta$-deformation and $x'=qx+\eta$ for the
$(q,\eta)$-deformation.

\textit{Remark}. There is an explicite formula for the elliptic
deformation of $\partial_x$.

\subsection{Standard deformations of functions}
Let $\exp(x)$ be the usual exponential of $x$:
\begin{equation}
\exp(x)=1+x+\frac{1}{2!}x^2 +\cdots+\frac{1}{n!}x^n+\cdots~.
\label{tolsDF13}
\end{equation}
and let $F_{n,m}(x)$ be a standard hypergeometric series:
\begin{equation}
F_{n,m}(x)=\sum_{k=0}^{\infty}
\frac{(a_1)_{(k)}(a_2)_{(k)}\cdots(a_n)_{(k)}}
{k!(b_1)_{(k)}(b_2)_{(k)}\cdots(b_m)_{(k)}}x^k~,\label{tolsDF14}
\end{equation}
where
\begin{equation}
(a)_{(k)}=a(a+1)\cdots(a+k-1)~.\label{tolsDF15}
\end{equation}

a). The trigonometric deformation of $\exp(x)$ and $F_{n,m}(x)$.
\begin{equation}
\exp_q(x)=1+x+\frac{1}{(2)_q!}x^2
+\cdots+\frac{1}{(n)_q!}x^n+\cdots~, \label{tolsDF16}
\end{equation}
where
\begin{equation}
(n)_q!=(1)_q(2)_q\cdots(n)_q~. \label{tolsDF17}
\end{equation}

If we replace the parameters $a_i$, $b_j$, and $k!$ in
(\ref{tolsDF14}) by the $q$-analogs $(a_i)_q$, $(b_j)_q$ and
$(k)_q!$ we obtain the basic hypergeometric ($q$-hypergeometric)
series $F_{n,m}^{(q)}(z)$.

b). The rational deformation of $\exp(x)$ and $F_{n,m}(x)$.
\begin{equation}
\exp_{\eta}(x)=1+x+\frac{1}{(2)_{\eta}!}x^2
+\cdots+\frac{1}{(n)_{\eta}!}x^n+\cdots~, \label{tolsDF18}
\end{equation}
where
\begin{equation}
(n)_{\eta}!=(1)_{\eta}(2)_{\eta}\cdots(n)_{\eta}~.
\label{tolsDF19}
\end{equation}

If we replace the parameters $a_i$, $b_j$, and $k!$ in
(\ref{tolsDF14}) by the $\eta$-analogs $(a_i)_{\eta}$,
$(b_j)_{\eta}$ and $(k)_{\eta}!$ we obtain the basic
hypergeometric ($\eta$-hypergeometric) series
$F_{n,m}^{(\eta)}(z)$.

c). The rational-trigonometric deformation of $\exp(x)$ and
$F_{n,m}(x)$.
\begin{equation}
\exp_{q\eta}(x)=1+x+\frac{1}{(2)_{q\eta}!}x^2
+\cdots+\frac{1}{(n)_{q\eta}!}x^n+\cdots~, \label{tolsDF20}
\end{equation}
where
\begin{equation}
(n)_{q\eta}!=(1)_{q\eta}(2)_{q\eta}\cdots(n)_{q\eta}~.
\label{tolsDF21}
\end{equation}

The replacement of $a_i$, $b_j$, and $k!$ in (\ref{tolsDF14}) by
the $(q,\eta)$-analogs $(a_i)_{q\eta}$, $(b_j)_{q\eta}$ and
$(k)_{q\eta}!$ gives us the $(q,\eta)$-hypergeometric series
\begin{equation}
F_{n,m}^{(q,\eta)}(x)=\sum_{k=0}^{\infty}
\frac{((a_1)_{q\eta})_{(k)}((a_2)_{q\eta})_{(k)}\cdots
((a_n)_{q\eta})_{(k)}}{(k)_{q\eta}!((b_1)_{q\eta})_{(k)}
((b_2)_{q\eta})_{(k)}\cdots((b_m)_{q\eta})_{(k)}}x^k~,
\label{tolsDF22}
\end{equation}
where
\begin{equation}
((a)_{q\eta})_{(k)}=(a)_{q\eta}(a+1)_{q\eta}
\cdots(a+k-1)_{q\eta}~.\qquad{} \label{tolsDF23}
\end{equation}
Setting here $q=1$ we obtain the $\eta$-hypergeometric series
$F_{n,m}^{(\eta)}(x)$. We can also introduce the $\eta$- and
($q,\eta$)- analogs of other special functions.

All deformed exponentials $\exp_q(x)$, $\exp_{\eta}(x)$ and
$\exp_{q\eta}(x)$ can be obtained from the functional equation
$f_1(x+y)=f_2(y)f_3(x)$ with regular functions $f_i(z)$ satisfying
the initial conditions $f_i(0)=1$, where the variables $x$ and $y$
satisfy the relations (\ref{tolsDF5})--(\ref{tolsDF7}) (see
\cite{tolsT2}).

\textit{Remark}. Analogous formulas exist also for the elliptic
case they can be found in \cite{tolsSp}.

\subsection{Standard deformations of universal enveloping algebras}
Let $g$ be a finite-dimensional Lie algebra, $g[u,u^{-1}]$ be a
loop algebra and $g[u]$ be a non-negative loop algebra over $g$.
We denote by $U(g)$, $U(g[u])$ and $U(g[u,u^{-1}])$ their
universal enveloping algebras.

a). The trigonometric deformation ($q$-deformation) of $U(g)$,
$U(g[u])$ and $U(g[u,u^{-1}])$ are well-known. They are denoted by
$U_q(g)$, $U_q(g[u])$ and $U_q(g[u,u^{-1}])$ \cite{tolsD}.

b). The rational deformation of $U_q(g[u])$ and $U_q(g[u,u^{-1}])$
are well-known. They are Yangian $Y_{\eta}(g)$ \cite{tolsD} and
its double $DY_{\eta}(g)$ \cite{tolsKT}.

c). The rational-trigonometric deformation of $U_q(g[u])$ and
$U_q(g[u,u^{-1}])$ are Drinfeldian $D_{q\eta}(g)$ \cite{tolsT1},
\cite{tolsT2} and its double $DD_{q\eta}(g)$.

\textit{Remark}. An elliptic deformation exists only for
$U(sl_n[u])$. It is well-known too.

\section{Rational-trigonometric differential  Knizhnik-Zamo\-lodchikov
and dynamical equations}

We consider the case $g=gl_M$ although results of this section are
also valid for an arbitrary simple complex Lie algebra $g$. Let
$e_{ab}$, $a,b=1,\ldots, M$, be the standard Cartan-Weyl basis of
$gl_M$: $[e_{ab},e_{cd}]=\delta_{bc}e_{ad}-\delta_{ad}e_{cb}$. The
element $C_2:=\sum\limits_{a,b=1}^Ne_{ab}e_{ba}\in U(gl_M)$ is a
$gl_M$-scalar, i.e. $[C_2,x]=0$ for any $x\in gl_M$, and it is
called the second order Casimir element. The element
$\Omega:=\frac{1}{2}\Big(\Delta(C_2)-C_2\otimes{\rm id}-{\rm
id}\otimes C_2\Big)=\sum\limits_{a,b=1}^M e_{ab}\otimes
e_{ba}\subset U(gl_M)\otimes U(gl_M)$, where $\Delta$ is a trivial
co-product $\Delta(x)=x\otimes{\rm id}+{\rm id}\otimes x$
($\forall x\in gl_M$), is called the Casimir two-tensor. The
two-tensor can be presented in the form $\Omega=$ $\Omega^+
+\Omega^-$, where
$\Omega^+=\frac{1}{2}\sum\limits_{a}e_{aa}\otimes
e_{aa}+\sum\limits_{1\leq a<b\leq M}e_{ab}\otimes e_{ba}$ and
$\Omega^-=\frac{1}{2}\sum\limits_{a}e_{aa}\otimes
e_{aa}+\sum\limits_{1\leq a<b\leq M}e_{ba}\otimes e_{ab}$. Note
that $(\omega\otimes\omega)(\Omega^{\pm})=\Omega^{\mp}$, where
$\omega$ is the Cartan automorphism: $\omega(e_{ab})=-e_{ba}$.

For any $x\in U(gl_M)$ we set $x_{(i)}=\underbrace{\mathop{{\rm
id}\otimes\cdots\otimes\textrm{id}\otimes\mathop{x}
\otimes\,\textrm{id}\otimes\cdots\otimes\textrm{id}}}_{N-times}^{\;\,i-th}$.
We consider $U(gl_M)$ as a subalgebra of $(U(gl_M))^{\otimes N}$,
the embedding $U(gl_M)\hookrightarrow (U(gl_M))^{\otimes N}$ being
given by the $N$-fold co-product, that is $x\mapsto
\Delta^N(x)=x_{(1)}+\ldots +x_{(N)}$ for any $x\in gl_M$.

For a nonzero complex number $\kappa$ we consider differential
operators $\nabla^{(r)}_{z_1^{}},\ldots,\nabla^{(r)}_{z_N^{}}$ and
$\nabla^{(t)}_{z_1^{}}, \ldots,\nabla^{(t)}_{z_N^{}}$ with
coefficient in $(U(gl_M))^{\otimes N}$ depending on complex
variable $z_1^{},\ldots,z_N^{}$ and
$\lambda_1^{},\cdots,\lambda_M^{}$ (see \cite{tolsTV1} and
\cite{tolsTV2}):
\begin{equation}
\nabla^{(r)}_{z_i}(z;\lambda)=\kappa\frac{\partial}{\partial
z_i}-\sum_{a=1}^{M}\lambda_a(e_{aa})_{(i)}-\sum_{j=1\atop j\neq
i}^{N}\frac{\Omega_{(ij)}}{z_i-z_j}~, \label{tolsDF27}
\end{equation}
\begin{eqnarray}
\nabla^{(t)}_{z_i}(z;\lambda)\!\!&=&\!\!\kappa z_i\frac{\partial}{\partial
z_i}-\sum_{a=1}^{M}(\lambda_a-e_{aa})(e_{aa})_{(i)}
-\sum_{j=1\atop i\neq i}^{N}\frac{z_i\Omega^{+}_{(ij)}+
z_j\Omega^{-}_{(ij)}} {z_i-z_j}~. \label{tolsDF28}
\end{eqnarray}
The operators $\nabla^{(r)}_{z_1^{}},\ldots,\nabla^{(r)}_{z_N^{}}$
(resp.$\nabla^{(t)}_{z_1^{}}, \ldots,\nabla^{(t)}_{z_N^{}}$) are
called the rational (resp. trigonometric) Knizhnik-Zamolodchikov
(KZ) operators. We set
\begin{equation}
\nabla^{(rt)}_{z_i}(z;\lambda)=\hbar\nabla^{(t)}_{z_i+
\frac{\eta}{\hbar}} (z+\mbox{\large$\frac{\eta}{\hbar}$};
(1+\mbox{\large$\frac{\eta}{\hbar}$})\lambda)~. \label{tolsDF29}
\end{equation}
It is easy to see that
\begin{equation}
\nabla^{(rt)}_{z_i}(z;\lambda)=\hbar\nabla^{(t)}_{z_i}(z;\lambda)+
\eta\nabla^{(r)}_{z_i}(z;\lambda)~. \label{tolsDF30}
\end{equation}
The rational (trigonometric or rational-trigonometric) KZ
equations is the system of the differential equations
\begin{equation}
\nabla^{(\alpha)}_{z_i}(z;\lambda)u(z;\lambda)=0~,\qquad
i=1,\ldots,N;\;\;{\alpha}=(r)\,((t)\;{\rm or}\;(rt))
\label{tolsDF31}
\end{equation}
for a function $u(z;\lambda):=u(z_{1}^{},\ldots,z_{N}^{};
\lambda_{1}^{},\cdots,\lambda_{M}^{})$ taking values in an
$N$-fold tensor product of $gl_M$-modules.

We also consider the differential operators
$D^{(r)}_{\lambda_{1}^{}},\ldots,D^{(r)}_{\lambda_{M}^{}}$ and
$D^{(t)}_{\lambda_{1}^{}},\ldots,D^{(t)}_{\lambda_{M}^{}}$ with
coefficient in $(U(gl_M))^{\otimes N}$ depending on complex
variables $z_1^{},\ldots,z_N^{}$ and
$\lambda_1^{},\cdots,\lambda_M^{}$ \cite{tolsTV1}:
\begin{equation}
D^{(r)}_{\lambda_a}(z;\lambda)=\kappa\frac{\partial}{\partial
\lambda_a}-\sum_{i=1}^{N}z_i(e_{aa})_{(i)}-\sum_{b=1\atop b\neq
a}^{M}\frac{e_{ab}e_{ba}-e_{aa}}{\lambda_a-\lambda_b}~,\phantom{aaaaaaaaa}
\label{tolsDF32}
\end{equation}
\begin{eqnarray}
D^{(t)}_{\lambda_a}(z;\lambda)\!\!&=&\!\!
\kappa\lambda_a\frac{\partial}{\partial
\lambda_a}+\frac{a_{aa}^2}{2}-\sum\limits_{i=1}^{N}z_i(e_{aa})_i-
\nonumber
\\
&& -\sum\limits_{b=1}^M\sum\limits_{1\leq
i<j\leq N}(e_{ab})_{(i)}(e_{ba})_{(j)}- \sum\limits_{b=1\atop
b\neq a}^{M}\frac{\lambda_b(e_{ab}e_{ba}-e_{aa})}
{\lambda_a-\lambda_b}~. \label{tolsDF33}
\end{eqnarray}
Remind that $e_{ab}=\sum\limits_{i=1}^{M}(e_{ab})_{(i)}$. The
operators $D^{(r)}_{z_1}, \ldots D^{(r)}_{z_M}$ (resp.
$D^{(t)}_{z_1}, \ldots D^{(t)}_{z_M}$) are called the rational
(resp. trigonometric) differential dynamic (DD) operators. We set
\begin{equation}
D^{(rt)}_{\lambda_a}(z;\lambda):=\hbar
D^{(t)}_{\lambda_a+\frac{\eta}{\hbar}}
((1+\mbox{\large$\frac{\eta}{\hbar}$})z;
\lambda+\mbox{\large$\frac{\eta}{\hbar}$})~. \label{tolsDF34}
\end{equation}
It is easy to see that
\begin{equation}
D^{(rt)}_{z_i}(z;\lambda)=\hbar D^{(t)}_{\lambda_a}(z;\lambda)+
\eta D^{(r)}_{\lambda_a}(z;\lambda) ~. \label{tolsDF35}
\end{equation}
The rational (trigonometric or rational-trigonometric) DD
equations is the system of the differential equations
\begin{equation}
D^{(\alpha)}_{\lambda_a}(z;\lambda)u(z;\lambda)=0~,\qquad
a=1,\ldots,M;\;\;{\alpha}=(r)\,((t)\;{\rm
or}\;(rt))\label{tolsDF36}
\end{equation}
for a function $u(z;\lambda):=u(z_{1}^{},\ldots,z_{N}^{};
\lambda_{1}^{},\cdots,\lambda_{M}^{})$ taking values in an
$N$-fold tensor product of $gl_M$-modules.

>From (\ref{tolsDF30}), (\ref{tolsDF35}) and Theorem 5.8 of the
paper \cite{tolsTV1} for any $i=1,\ldots, N$ and $a=1,\ldots, M$
we have the duality:
\begin{equation}
\nabla^{(rt)}_{z_i}(z;\lambda)_{M}^{}\simeq
D^{(rt)}_{z_i}(\lambda;z)_{N}^{}~,\qquad
D^{(rt)}_{\lambda_a}(z;\lambda)_{M}^{}\simeq
\nabla^{(rt)}_{\lambda_a}(\lambda;z)_{N}^{}~,\label{tolsDF37}
\end{equation}
where the index $M$ ($N$) is connected with the Lie algebra
$gl_M$($gl_N$).

\vskip 7pt \noindent \textit{Acknowledgment}. This work was
supported by Russian Foundation for Fundamental Research, grant
No. RFBR-02-01-00668, and CRDF RMI-2334-MO-02.


\end{document}